# OPTIMAL TESTING OF EQUIVALENCE HYPOTHESES

By Joseph P. Romano

*Stanford University*

In this paper we consider the construction of optimal tests of equivalence hypotheses. Specifically, assume $X_1, \ldots, X_n$ are i.i.d. with distribution $P_\theta$, with $\theta \in \mathbb{R}^k$. Let $g(\theta)$ be some real-valued parameter of interest. The null hypothesis asserts $g(\theta) \notin (a,b)$ versus the alternative $g(\theta) \in (a,b)$. For example, such hypotheses occur in bioequivalence studies where one may wish to show two drugs, a brand name and a proposed generic version, have the same therapeutic effect. Little optimal theory is available for such testing problems, and it is the purpose of this paper to provide an asymptotic optimality theory. Thus, we provide asymptotic upper bounds for what is achievable, as well as asymptotically uniformly most powerful test constructions that attain the bounds. The asymptotic theory is based on Le Cam's notion of asymptotically normal experiments. In order to approximate a general problem by a limiting normal problem, a UMP equivalence test is obtained for testing the mean of a multivariate normal mean.

**1. Introduction.** Suppose $X_1, \ldots, X_n$ are i.i.d. with distribution $P_\theta$, where $\theta$ is a vector in $\mathbb{R}^k$. Let $g(\theta)$ be some real-valued parameter of interest. The hypothesis testing problem we study in this paper is of the following form: the null hypothesis asserts $g(\theta) \notin (a,b)$ versus the alternative $g(\theta) \in (a,b)$.

This setup arises when trying to demonstrate equivalence (or sometimes called bioequivalence) of treatments. By comparing a pharmacokinetic parameter of a new drug to the standard drug, bioequivalence is declared if the parameter $g(\theta)$ lies in the interval $(a,b)$, where $(a,b)$ is specified by a regulatory agency. For example, if $g(\theta)$ is the difference in treatment means, then equivalence corresponds to values of $g(\theta)$ near zero, and so $(a,b) = (-\Delta, \Delta)$ for some $\Delta > 0$. Then, rejection of the null hypothesis is the same as declaring equivalence. By formulating the null hypothesis as $g(\theta) \notin (-\Delta, \Delta)$, the risk of marketing an alternative drug that does not behave like the standard drug is controlled. In some situations it may be more appropriate to specify equivalence by a ratio of means, and equivalence is then declared if the ratio is near one, so that $(a, b)$ would be an interval containing one. More generally, the problem may consist of determining equivalence across several parameters, but only the simple real-valued case is treated here. A very nice recent account of testing hypotheses of equivalence is given in [17]. For the remainder of the paper, we assume without loss of generality that $(a, b) = (-\Delta, \Delta)$.

If $P$ belongs to a one-parameter exponential family, then a (uniformly most powerful) UMP level $\alpha$ test exists; see [10], Theorem 6 in Chapter 3. More generally, if the family of distributions is strictly positive of order 3 and other mild continuity conditions are satisfied, then a UMP test exists; the more general result appears as Problem 30 in Chapter 3 of [10] and is proved in [7].

EXAMPLE 1.1 (*Normal location model*). First, suppose $X$ is $N(\mu, \sigma^2)$, the normal distribution with mean $\mu$ and variance $\sigma^2$. Assume $\sigma^2$ is known. The problem is to test $|\mu| \geq \Delta$ versus $|\mu| < \Delta$. Applying the previously mentioned result, the UMP level $\alpha$ test rejects if $|X| < C$, where $C = C(\alpha, \Delta, \sigma)$ satisfies

$$\Phi\left(\frac{C - \Delta}{\sigma}\right) - \Phi\left(\frac{-C - \Delta}{\sigma}\right) = \alpha \tag{1}$$

and $\Phi(\cdot)$ is the standard normal c.d.f.

Next, suppose $X_1, \ldots, X_n$ are i.i.d. $N(\mu, \sigma^2)$. For testing the same hypothesis, the UMP level $\alpha$ test rejects if $n^{1/2}|\bar{X}_n| \leq C(\alpha, n^{1/2}\Delta, \sigma)$.

If $\sigma$ is unknown, no UMP test exists, nor do unbiasedness or invariance considerations lead to an optimal test.

Outside a small class of models, no optimality theory is available for tests of equivalence. Wellek [17] provides a general construction of asymptotically valid tests, based on some asymptotically normal estimators, but no theory is provided to prove optimality of such procedures. The main goal of this paper is to provide an asymptotic optimality theory for such problems. Specifically, we obtain bounds for the asymptotic power of tests of equivalence for a large class of models, as well as construct efficient tests that attain these bounds. As will be seen, the results flow from Le Cam's approach based on convergence of experiments; see [8]. In order for this approach to be viable, we need to determine an optimal test for the limiting normal experiment; this is accomplished in Section 2, where an exact finite sample theory uniformly most powerful test is derived for testing the equivalence of a linear function of a multivariate normal mean.

In Section 3 we consider asymptotic efficiency. We will formulate the asymptotic problem in two distinct ways. First, consider the case when the null hypothesis parameter space is the complement of a fixed interval $(-\Delta, \Delta)$ is considered. Then, we analyze the case when this interval changes (and shrinks) with $n$. In each case, attainable upper bounds for the asymptotic power of tests are obtained. The upper bounds in the two approaches actually differ, and we prefer the second approach.

In fact, Janssen [6] has already considered the problem of testing equivalence in a semiparametric two-sample framework, which in many ways is a more difficult problem. He too considered a shrinking alternative parameter space. Building on the work of Pfanzagl [12, 13] and Janssen [5], his technique also relies on a reduction to an asymptotically normal experiment. However, he proves optimality of the power function at a particular value of the functional of interest, which, in our case, corresponds to the value of $g(\theta)$ being zero. He also imposes an asymptotic similarity condition (see (36) of [6]). Here, we obtain asymptotically optimal tests, uniformly over values of the parameters in the region for which equivalence is declared, and no asymptotic unbiasedness condition is imposed. Of course, we are working in a parametric framework in order to obtain such a clean result. But, our results can be used in semiparametric models by an appropriate reduction to a parametric least favorable submodel. Some asymptotic nonparametric results which differ slightly from those of Janssen [6] can be obtained from the author.

Even in the normal one-sample problem with unknown variance, the problem of testing for equivalence has a rich history, and some of the literature is given in Example 3.1. Our asymptotic results will apply quite generally to parametric models, under the weak assumption of quadratic mean differentiability. The results generalize immediately to two-sample (or $s$-sample) problems, as well as to more complicated designs (such as a crossover design), as long as the underlying model is smooth enough to permit convergence to a normal experiment, since the optimal test in the limiting normal experiment is completely specified in Theorem 2.1.

**2. A finite sample UMP test.** Throughout, $\Phi(\cdot)$ is the standard normal distribution function and $z_\alpha$ satisfies $\Phi(z_\alpha) = \alpha$. Before discussing optimality for general models, we first need to derive the optimal test in the appropriate limiting normal experiment. This is obtained in the following result.

THEOREM 2.1. *Suppose $(X_1, \ldots, X_k)$ is multivariate normal $N(\mu, \Sigma)$ with unknown mean $\mu = (\mu_1, \ldots, \mu_k)^T$ and known covariance matrix $\Sigma$ (possibly nonsingular). Fix $\delta > 0$ and any vector $a = (a_1, \ldots, a_k)^T$ satisfying $a^T \Sigma a > 0$.*

*Consider testing*

$$H : \left|\sum_{i=1}^{k} a_i \mu_i\right| \geq \delta \quad vs. \quad K : \left|\sum_{i=1}^{k} a_i \mu_i\right| < \delta.$$

*Then a UMP level $\alpha$ test exists and it rejects $H$ if*

$$\left|\sum_{i=1}^{k} a_i X_i\right| < C,$$

*where $C = C(\alpha, \delta, \sigma)$ satisfies*

(2) $$\Phi\left(\frac{C-\delta}{\sigma}\right) - \Phi\left(\frac{-C-\delta}{\sigma}\right) = \alpha$$

*and $\sigma^2 = a^T \Sigma a$. Hence, the power of this test against an alternative $(\mu_1, \ldots, \mu_k)$ with $|\sum_i a_i \mu_i| = \delta' < \delta$ is*

$$\Phi\left(\frac{C-\delta'}{\sigma}\right) - \Phi\left(\frac{-C-\delta'}{\sigma}\right).$$

PROOF. The proof will consider four cases in increasing generality.

*Case* 1. Suppose $k = 1$, so that $X_1 = X$ is $N(\mu, \sigma^2)$ and we are testing $|\mu| \geq \delta$ versus $|\mu| < \delta$. Fix an alternative $\mu = m$ with $|m| < \delta$. Reduce the composite null hypothesis to a simple one via a least favorable distribution that places mass $p$ on $N(\delta, \sigma^2)$ and mass $1 - p$ on $N(-\delta, \sigma^2)$. The value of $p$ will be chosen shortly so that such a distribution is least favorable (and will be seen to depend on $m$, $\alpha$, $\sigma$ and $\delta$). By the Neyman–Pearson lemma, the MP test of

$$pN(\delta, \sigma^2) + (1-p)N(-\delta, \sigma^2) \quad vs. \quad N(m, \sigma^2)$$

rejects for small values of

(3) $$\frac{p \exp[-(1/(2\sigma^2))(X-\delta)^2] + (1-p)\exp[-(1/(2\sigma^2))(X+\delta)^2]}{\exp[-(1/(2\sigma^2))(X-m)^2]},$$

or, equivalently, for small values of $f(X)$, where

$$f(x) = p \exp[(\delta - m)X/\sigma^2] + (1-p)\exp[-(\delta + m)X/\sigma^2].$$

We can now choose $p$ so that $f(C) = f(-C)$, so that $p$ must satisfy

(4) $$\frac{p}{1-p} = \frac{\exp[(\delta+m)C/\sigma^2] - \exp[-(\delta+m)C/\sigma^2]}{\exp[(\delta-m)C/\sigma^2] - \exp[-(\delta-m)C/\sigma^2]}.$$

Since $\delta - m > 0$ and $\delta + m > 0$, both the numerator and denominator of the right-hand side of (4) are positive, so the right-hand side is a positive number; but, $p/(1-p)$ is a nondecreasing function of $p$ with range $[0,\infty)$ as $p$ varies from 0 to 1. Thus, $p$ is well defined. Also, observe $f''(x) \geq 0$ for all $x$. It follows that (for this special choice of $C$)

$$\{X : f(X) \leq f(C)\} = \{X : |X| \leq C\}$$

is the rejection region of the MP test. Such a test is easily seen to be level $\alpha$ for the original composite null hypothesis because its power function is symmetric and decreases away from zero. Thus, the result follows by Theorem 6 in Section 3.7 of [10].

*Case* 2. Consider now general $k$, so that $(X_1, \ldots, X_k)$ has mean $(\mu_1, \ldots, \mu_k)$ and covariance matrix $\Sigma$. However, consider the special case $(a_1, \ldots, a_k) = (1, 0, \ldots, 0)$, so we are testing $|\mu_1| \geq \delta$ versus $|\mu_1| < \delta$. Also, assume $X_1$ and $(X_2, \ldots, X_k)$ are independent, so that the first row and first column of $\Sigma$ are zero except the first entry, which is $\sigma^2$ (assumed positive). Using the same reasoning as in Case 1, fix an alternative $m = (m_1, \ldots, m_k)$ with $|m_1| < \delta$ and consider testing

$$pN((\delta, m_2, \ldots, m_k), \Sigma) + (1-p)N((-\delta, m_2, \ldots, m_k), \Sigma)$$

versus $N((m_1, \ldots, m_k), \Sigma)$. The likelihood ratio is, in fact, the same as (3) because each term is now multiplied by the density of $(X_2, \ldots, X_k)$ (by independence), and these densities cancel. The UMP test from Case 1, which rejects when $|X_1| \leq C$, is UMP in this situation as well.

*Case* 3. As in Case 2, consider $a_1 = 1$ and $a_i = 0$ if $i > 1$, but now allow $\Sigma$ to be arbitrary. Reduce the problem to Case 2 by an appropriate linear transformation. Simply let $Y_1 = X_1$ and, for $i > 1$, let

$$Y_i = X_i - \frac{\text{Cov}(X_1, X_i)}{\text{Var}(X_1)} X_1,$$

so that $\text{Cov}(Y_1, Y_i) = 0$ if $i > 0$. Thus, the problem of testing $E(Y_1) = E(X_1)$, based on $Y = (Y_1, \ldots, Y_k)$, is in the form studied in Case 2, and the UMP test rejects for small $|Y_1| = |X_1|$.

*Case* 4. Now, consider arbitrary $(a_1, \ldots, a_k)$ satisfying $a^T \Sigma a > 0$. Let $Z = OX$, where $O$ is any orthogonal matrix with first row $(a_1, \ldots, a_k)$. Then $E(Z_1) = \sum_{i=1}^k a_i \mu_i$, and the problem of testing $|E(Z_1)| \geq \delta$ versus $|E(Z_1)| < \delta$ reduces to Case 3. Hence, the UMP test rejects for small values of $|Z_1| = |\sum_{i=1}^k a_i X_i|$. □

Next, we summarize some simple but useful properties of the critical constants $C(\alpha, \delta, \sigma)$ and the optimal power of the above UMP test that will be used later.

REMARK 2.1. It is easy to check that, as a function of $C$, the function $h(C)$ given by

$$h(C) = \Phi\left(\frac{C-\delta}{\sigma}\right) - \Phi\left(\frac{-C-\delta}{\sigma}\right)$$

is increasing in $C$. Since

$$h(\delta - \sigma z_{1-\alpha}) = \alpha - \Phi(-2\sigma^{-1}\delta + z_{1-\alpha}) < \alpha,$$

it follows that

(5) $$C(\alpha, \delta, \sigma) > \delta - \sigma z_{1-\alpha}.$$

REMARK 2.2. The function $C(\alpha, \delta, \sigma)$ satisfies

(6) $$\frac{C(\alpha, \delta, \sigma)}{\sigma} = C\left(\alpha, \frac{\delta}{\sigma}, 1\right).$$

It is also easy to check that

$$C(\alpha, \varepsilon, 1) \to z_{(1-\alpha/2)}$$

as $\varepsilon \to 0$ and $C(\alpha, B, 1) \to \infty$ as $B \to \infty$.

REMARK 2.3. For fixed $C = C(\alpha, \delta, \sigma)$, the function

$$f(\gamma) = \Phi(C - \gamma) - \Phi(-C - \gamma)$$

is decreasing in $\gamma$; to see why, just differentiate $f$. So if $0 \leq \gamma < \delta$, then $f(\gamma) > \alpha$.

**3. Asymptotic optimality.** Throughout this section, we assume $X_1, \ldots, X_n$ are i.i.d. according to a distribution $P_\theta$, with $\theta \in \Omega$, where $\Omega$ is an open subset of $\mathbb{R}^k$. Here the observations take values in a sample space $\mathcal{X}$. Assume $P_\theta$ has density $p_\theta$ with respect to $\mu$. We will assume the family is quadratic mean differentiable (q.m.d.) at certain values $\theta_0$ of $\theta$; that is, there exists a vector of real-valued functions $\eta(\cdot, \theta_0) = (\eta_1(\cdot, \theta_0), \ldots, \eta_k(\cdot, \theta_0))^T$ such that

(7) $$\int_{\mathcal{X}} [\sqrt{p_{\theta_0+h}(x)} - \sqrt{p_{\theta_0}(x)} - \langle \eta(x, \theta_0), h \rangle]^2 \, d\mu(x) = o(|h|^2)$$

as $|h| \to 0$. Here $h$ is a vector in $\mathbb{R}^k$ and $|h|$ denotes its Euclidean norm. For such a family the Fisher information matrix at $\theta_0$ is the matrix $I(\theta_0)$ with $(i, j)$ entry

$$I_{i,j}(\theta_0) = 4 \int \eta_i(x, \theta_0) \eta_j(x, \theta_0) \, d\mu(x).$$

We also define the *score* vector $Z_n$ to be

$$(8) \quad Z_n = Z_n(\theta_0) = 2n^{-1/2} \sum_{i=1}^{n} [\eta(X_i, \theta_0) / p_{\theta_0}^{1/2}(X_i)].$$

For a review of families that are q.m.d., as well as the history and importance of this notion, see [9]. In particular, we make heavy use of the fact that such families are locally asymptotically normal, and so the testing problem under consideration can be approximated by a certain normal testing problem.

Interest focuses on $g(\theta)$, where $g$ is a function from $\Omega$ to $\mathbb{R}$. Assume $g$ is differentiable with gradient vector $\dot{g}(\theta)$ of dimension $1 \times k$. We will formulate the problem in two distinct ways. First, we consider the case when the null hypothesis parameter space $\Omega_0$ is the complement of a fixed interval $(-\Delta, \Delta)$. Then, we study the case when this interval changes (and shrinks) with $n$.

3.1. *Fixed parameter spaces.* Fix $\Delta > 0$. The problem is to test $|g(\theta)| \geq \Delta$ versus $|g(\theta)| < \Delta$. We implicitly assume $g$ is such that there exists a $\theta$ such that $g(\theta) \geq \Delta$, as well as a $\theta$ with $g(\theta) \leq -\Delta$. For any fixed alternative value $\theta$ with $|g(\theta)| < \Delta$, the power of any reasonable test against $\theta$ will tend to one. Therefore, as is customary (see [16], Chapters 14 and 15), we compare power functions at local alternatives. Consider any fixed $\theta_0$ satisfying $|g(\theta_0)| = \Delta$. For sake of argument, consider the case $g(\theta_0) = -\Delta$. In order to derive an (obtainable) upper bound for the limiting power of a test sequence $\phi_n$ under $\theta_0 + hn^{-1/2}$, a crude way to bound the power is based on the simple fact that any level $\alpha$ test for testing $|g(\theta)| \geq \Delta$ versus $|g(\theta)| < \Delta$ is also level $\alpha$ for testing $g(\theta) \leq -\Delta$ versus $g(\theta) > -\Delta$. Since upper bounds for the (asymptotic) power are well known for the latter testing problem (as in [16], Theorem 15.4), an immediate result follows. In this asymptotic setup, the statistical problem is somewhat degenerate, as it becomes one of testing a one-sided hypothesis. For example, suppose $X_1, \ldots, X_n$ are i.i.d. $N(\theta, 1)$. Then for large $n$, one can distinguish $\theta \leq -\Delta$ and $\theta > -\Delta$ with error probabilities that are uniformly small and tend to zero exponentially fast with $n$. In essence, the statistical issue arises only if the true $\theta$ is near the boundary of $[-\Delta, \Delta]$, in which case determining significance essentially becomes one of testing a one-sided hypothesis.

EXAMPLE 3.1 (*Normal one-sample problem*). Suppose $X_1, \ldots, X_n$ are i.i.d. $N(\mu, \sigma^2)$, with both parameters unknown. Consider testing $|\mu| \geq \Delta$ versus $|\mu| < \Delta$. The standard $t$-test for testing the one-sided hypothesis $\mu \leq -\Delta$ against $\mu > -\Delta$ rejects if

$$n^{1/2}(\bar{X}_n + \Delta)/S_n > t_{n-1,1-\alpha},$$

where $S_n^2$ is the (unbiased) sample variance and $t_{n-1,1-\alpha}$ is the $1-\alpha$ quantile of the $t$-distribution with $n-1$ degrees of freedom. Similarly, the standard $t$-test of the hypothesis $\mu \geq \Delta$ rejects if

$$n^{1/2}(\bar{X}_n - \Delta)/S_n < -t_{n-1,1-\alpha}.$$

The intersection of these rejection regions is therefore a level $\alpha$ test of the null hypothesis $|\mu| \geq \Delta$. Such a construction that intersects the rejection regions of two one-sided tests (TOST) was proposed in [18] and [15], and can be seen as a special case of Berger's [2] intersection-union tests; see [3] for a review. The resulting TOST is given by the test $\phi_n^{\text{TOST}}$ that rejects when $|\bar{X}_n| < \Delta - n^{-1/2} S_n t_{n-1,1-\alpha}$. The asymptotic power of $\phi_n^{\text{TOST}}$ against a sequence with mean $-\Delta + hn^{-1/2}$ ($h > 0$) and variance fixed at $\sigma^2$ can be calculated directly as

$$P_{\Delta + hn^{-1/2},\sigma}\{|\bar{X}_n| < \Delta - n^{-1/2} S_n t_{n-1,1-\alpha}\} = \Phi\left(z_{1-\alpha} - \frac{h}{\sigma}\right),$$

which is the optimal bound for the one-sided testing problem given in Theorem 15.4 of [16]. A similar calculation applies to sequences of the form $\Delta - hn^{-1/2}$. Thus, the TOST is asymptotically optimal in this setup. It should be remarked that the TOST has been criticized because it is biased (in finite samples) and tests have been proposed that have greater power; some proposals are reviewed and studied in [3], [4] and [11]. These points are valid, but no test can have greater asymptotic power against such local alternatives. On the other hand, the TOST will be seen to be inefficient under the asymptotic setup of the next section.

3.2. *Shrinking alternative parameter space.* We now consider a second asymptotic formulation of the problem. The null hypothesis asserts $|g(\theta)| \geq \delta n^{-1/2}$ and the alternative hypothesis asserts $|g(\theta)| < \delta n^{-1/2}$. Notice now that, in this asymptotic study, the parameter spaces (or hypotheses) are changing with $n$. Of course, a given hypothesis testing situation deals with a particular $n$, and there is flexibility in how the problem is embedded into a sequence of similar problems to get a useful approximation. Indeed, if equivalence corresponds to $|g(\theta)| < \Delta$, we can always set up the problem by choosing $\delta = \Delta n^{1/2}$. From an asymptotic point of view, it makes sense to allow the null hypothesis parameter space to change with $n$, or else the problem becomes degenerate in the sense that the values of $\Delta$ and $-\Delta$ for $g(\theta)$ can be perfectly distinguished asymptotically. In testing for bioequivalence, for example, $\Delta$ represents a small value so that a value of $|g(\theta)| \leq \Delta$ is deemed sufficiently close to zero in a clinical sense. In a particular situation (such as the previous example with $\sigma$ not too small), a value for $|g(\theta)|$ of $\Delta$ cannot be perfectly tested against a value of $g(\theta) = 0$. Thus, if 0 is in some sense not far from both $\Delta$ and $-\Delta$, it follows that $\Delta$ and $-\Delta$ are not far from each other either, and the asymptotic setup should reflect this.

We implicitly assume there exists some $\theta$ with $g(\theta) > 0$, as well as some $\theta$ with $g(\theta) < 0$. The main result of this section is the following theorem.

THEOREM 3.1. *Suppose $X_1, \ldots, X_n$ are i.i.d. according to $P_\theta$, $\theta \in \Omega$, where $\Omega$ is assumed to be an open subset of $\mathbb{R}^k$. Consider testing the null hypothesis*

$$\theta \in \Omega_{0,n} = \{\theta : |g(\theta)| \geq \delta n^{-1/2}\}$$

*versus $|g(\theta)| < \delta n^{-1/2}$, where the function $g$ from $\mathbb{R}^k$ to $\mathbb{R}$ is assumed differentiable with gradient $\dot{g}(\theta)$. Assume for every $\theta$ with $g(\theta) = 0$ that the family $\{P_\theta, \theta \in \Omega\}$ is q.m.d. at $\theta$ and $I(\theta)$ is nonsingular.*

(i) *Let $\phi_n = \phi_n(X_1, \ldots, X_n)$ be a uniformly asymptotically level $\alpha$ sequence of tests, so that*

$$\limsup_{n \to \infty} \sup_{\Omega_{0,n}} E_\theta(\phi_n) \leq \alpha.$$

*Assume $\theta_0$ satisfies $g(\theta_0) = 0$. Then, for any $h$ such that $|\langle \dot{g}(\theta_0)^T, h \rangle| = \delta' < \delta$,*

$$(9) \qquad \limsup_{n \to \infty} E_{\theta_0 + h n^{-1/2}}(\phi_n) \leq \Phi\left(\frac{C - \delta'}{\sigma_{\theta_0}}\right) - \Phi\left(\frac{-C - \delta'}{\sigma_{\theta_0}}\right),$$

*where $\sigma_{\theta_0}^2$ is given by*

$$(10) \qquad \sigma_{\theta_0}^2 = \dot{g}(\theta_0) I^{-1}(\theta_0) \dot{g}(\theta_0)^T$$

*and $C = C(\alpha, \delta, \sigma_{\theta_0})$ satisfies (2).*

(ii) *Let $\hat{\theta}_n$ be any estimator satisfying*

$$(11) \qquad n^{1/2}(\hat{\theta}_n - \theta_0) = I^{-1}(\theta_0) Z_n + o_{P_{\theta_0}^n}(1),$$

*(such as an efficient likelihood estimator). Suppose $I(\theta)$ is continuous in $\theta$ and $\dot{g}(\theta)$ is continuous at $\theta_0$. Then the test sequence $\phi_n$ that rejects when $n^{1/2}|g(\hat{\theta}_n)| \leq C(\alpha, \delta, \hat{\sigma}_n)$, where*

$$\hat{\sigma}_n^2 = \dot{g}(\hat{\theta}_n) I^{-1}(\hat{\theta}_n) \dot{g}(\hat{\theta}_n)^T,$$

*is pointwise asymptotically level $\alpha$ and is locally asymptotically UMP in the sense that the inequality (9) is an equality. In fact, the same properties hold for any test sequence that rejects if $|T_n| < C(\alpha, \delta, \hat{\sigma}_n)$, if $T_n$ satisfies*

$$T_n = \dot{g}(\theta_0) I^{-1}(\theta_0) Z_{n,\theta_0} + o_{P_{\theta_0}^n}(1)$$

*for every $\theta_0 \in \Omega_0$, where $Z_{n,\theta_0}$ is the score vector defined in (8).*

PROOF. Fix $\theta_0$ satisfying $g(\theta_0) = 0$. We will derive an upper bound for the limiting power of a test sequence $\phi_n$ under $\theta_0 + hn^{-1/2}$. Note that

$$g(\theta_0 + hn^{-1/2}) = n^{-1/2}\langle \dot{g}(\theta_0)^T, h\rangle + o(n^{-1/2}).$$

So, if $h$ is such that $|\langle \dot{g}(\theta_0)^T, h\rangle| > \delta$, then $|g(\theta_0 + hn^{-1/2})| > \delta n^{-1/2}$ for all sufficiently large $n$. Hence, if $\phi_n$ has limiting size $\alpha$, then, for such an $h$,

(12) $$\limsup_{n \to \infty} E_{\theta_0 + hn^{-1/2}}(\phi_n) \leq \alpha.$$

Since the family is q.m.d., the sequence of experiments $P^n_{\theta_0+hn^{-1/2}}$ (indexed by a vector $h$) converges to a limiting (multivariate) normal experiment with unknown mean vector $h$ and known covariance matrix $I^{-1}(\theta_0)$. Therefore, we can approximate the power of a test sequence $\phi_n$ by the power of a test $\phi = \phi(X)$ for the (limit) experiment based on $X$ from the model $N(h, I^{-1}(\theta_0))$; see Lemma 3.4.4 of [14] or Theorem 15.1 of [16]. So, let $\beta_\phi(h)$ denote the power function of $\phi(X)$ when $X \sim N(h, I^{-1}(\theta_0))$. Then (12) implies $\beta_\phi(h) \leq \alpha$ if $|\langle \dot{g}(\theta_0)^T, h\rangle| > \delta$. By continuity of $\beta_\phi(h)$, $\beta_\phi(h) \leq \alpha$ for any $h$ with $|\langle \dot{g}(\theta_0)^T, h\rangle| \geq \delta$. The choice of $\phi$ to maximize $\beta_\phi(h)$ for this limiting normal problem was given in Theorem 2.1 with $\Sigma = I^{-1}(\theta_0)$ and $a^T = \dot{g}(\theta_0)$. Thus, if $\phi$ is level $\alpha$ for testing $|\langle \dot{g}(\theta_0)^T, h\rangle| \geq \delta$ and $h$ satisfies $|\langle \dot{g}(\theta_0)^T, h\rangle| = \delta' < \delta$, then

$$\beta_\phi(h) \leq \Phi\left(\frac{C - \delta'}{\sigma_{\theta_0}}\right) - \Phi\left(\frac{-C - \delta'}{\sigma_{\theta_0}}\right),$$

and $C = C(\alpha, \delta, \sigma_{\theta_0})$ satisfies (2).

To prove (ii), consider the test that rejects when $n^{1/2}|g(\hat{\theta}_n)| \leq C(\alpha, \delta, \hat{\sigma}_n)$. Fix $h$ such that $|\langle \dot{g}(\theta_0)^T, h\rangle| = \delta' < \delta$ and let $\theta_n = \theta_0 + hn^{-1/2}$. Then, using standard contiguity arguments, under $\theta_n$,

$$n^{1/2}[g(\hat{\theta}_n) - g(\theta_n)] \xrightarrow{N} (0, \sigma^2_{\theta_0}).$$

But

$$n^{1/2}g(\theta_n) = \langle h, \dot{g}(\theta_0)^T\rangle + o(1).$$

Therefore, under $\theta_n$,

$$n^{1/2}g(\hat{\theta}_n) \xrightarrow{L} N(\langle h, \dot{g}(\theta_0)^T\rangle, \sigma^2_{\theta_0}).$$

Also, under $\theta_n$, $\hat{\sigma}_n$ tends in probability to $\sigma_{\theta_0}$, and so $C(\alpha, \delta, \hat{\sigma}_n)$ tends in probability to $C(\alpha, \delta, \sigma_{\theta_0})$. Hence, letting $Z$ denote a standard normal variable,

$$P_{\theta_n}\{n^{1/2}|g(\hat{\theta}_n)| \leq C(\alpha, \delta, \hat{\sigma}_n)\} \to P\{|\sigma_{\theta_0}Z + \langle h, \dot{g}(\theta_0)^T\rangle| \leq C(\alpha, \delta, \sigma_{\theta_0})\},$$

which agrees with the right-hand side of (9). □

OPTIMAL EQUIVALENCE TEST    11EXAMPLE 3.2 (*Normal one-sample problem, Example* 3.1, *continued*).
Suppose $X_1,\ldots,X_n$ are i.i.d. $N(\mu,\sigma^2)$ with both parameters unknown, so that $\theta = (\mu,\sigma)$. Let $g(\theta) = \mu$ and consider testing $|\mu| \geq \delta n^{-1/2}$ versus $|\mu| < \delta n^{-1/2}$. By the previous theorem, for any test sequence $\phi_n$ with limiting size bounded by $\alpha$ and any $h$ with $|h| < \delta$,

$$\text{(13)} \qquad E_{hn^{-1/2},\sigma}(\phi_n) \leq \Phi\left(\frac{C-h}{\sigma}\right) - \Phi\left(\frac{-C-h}{\sigma}\right),$$

where $C = C(\alpha,\delta,\sigma)$ satisfies (2). A test whose limiting power achieves this bound is given by the test $\phi_n^*$ that rejects when

$$n^{1/2}|\bar{X}_n| \leq C(\alpha,\delta,S_n),$$

where $S_n^2$ is the (unbiased) sample variance (or any consistent estimator of $\sigma^2$). In the normal model, such an approximate test was first proposed by Anderson and Hauck [1] (but in a two-sample context); in essence, the general construction of Theorem 3.1 can be viewed as an extension of their method.

On the other hand, the test $\phi_n^{\text{TOST}}$ given in Example 3.1 is *no longer* asymptotically efficient. This test (with $\Delta = \delta n^{-1/2}$) rejects when

$$n^{1/2}|\bar{X}_n| < \delta - S_n t_{n-1,1-\alpha}$$

and has power against $(\mu,\sigma) = (hn^{-1/2},\sigma)$ given by

$$\text{(14)} \quad P_{hn^{-1/2},\sigma}\left\{\frac{-\delta + S_n t_{n-1,1-\alpha} - h}{\sigma} < Z_n < \frac{\delta - S_n t_{n-1,1-\alpha} - h}{\sigma}\right\},$$

where

$$Z_n = n^{1/2}(\bar{X}_n - hn^{-1/2})/\sigma \sim N(0,1).$$

Also, $S_n \to \sigma$ in probability and $t_{n-1,1-\alpha} \to z_{1-\alpha}$. By Slutsky's theorem, (14) converges to

$$\text{(15)} \qquad P\left\{\frac{-\delta}{\sigma} + z_{1-\alpha} - \frac{h}{\sigma} < Z < \frac{\delta}{\sigma} - z_{1-\alpha} - \frac{h}{\sigma}\right\},$$

where $Z \sim N(0,1)$. Observe that this last expression is positive only if $\sigma z_{1-\alpha} < \delta$; otherwise, the limiting power is zero! On the other hand, the limiting optimal power of $\phi_n^*$ is always positive (and greater than $\alpha$ when $|h| < \delta$). Even when the limiting power of $\phi_n^{\text{TOST}}$ is positive, it is always strictly less than that of $\phi_n^*$. Note that the limiting expression (15) for the power of $\phi_n^{\text{TOST}}$ corresponds exactly to using a TOST test in the limiting experiment $N(h,\sigma^2)$, where you are testing $|h| \geq \delta$ versus $|h| < \delta$ with $\sigma$ known; such a procedure is conservative and less powerful than the UMP test. In general, (5) implies that

$$C(\alpha,\delta,\hat{\sigma}_n) > \delta - \hat{\sigma}_n z_{1-\alpha},$$

which shows that the test $\phi_n^*$ of Theorem 3.1 is always more powerful than the asymptotic TOST construction.

REMARK 3.1. Thus far, we have considered testing for two situations, first where the null hypothesis specifies $|g(\theta)| \geq \Delta$ and next where $|g(\theta)| \geq \delta/n^{1/2}$. Of course, one can also consider the general situation where the null hypothesis is specified by $|g(\theta)| \geq \delta/\tau_n$, where $\tau_n \to \infty$. For the purposes of this discussion, suppose $g(\theta) = \theta \in \mathbb{R}$. So, suppose we are testing $|\theta| \geq \delta/\tau_n$ with $\tau_n \to \infty$ at a rate slower than $n^{1/2}$, so that $\tau_n = o(n^{1/2})$. By contiguity arguments, the optimal limiting power will be nondegenerate (meaning away from $\alpha$ and 1) for alternatives of the form $\delta/\tau_n - h/n^{1/2}$ or $\delta/\tau_n + h/n^{1/2}$ for $h > 0$, or, more generally, if $h/n^{1/2}$ is replaced by any sequence $\varepsilon_n$ satisfying $\varepsilon \asymp n^{-1/2}$. But, if $\tau_n = o(n^{1/2})$, then $\delta/\tau_n$ and $-\delta/\tau_n$ can be perfectly distinguished, and so we are essentially in the first asymptotic setup. That is, the asymptotically optimal power against an alternative sequence $\delta/\tau_n - h/n^{1/2}$ is the same as for testing a one-sided hypothesis $\theta \geq \delta/\tau_n$ versus $\theta < \delta/\tau_n$.

On the other hand, suppose $\tau_n \to \infty$ faster than $n^{1/2}$, so that $n^{1/2}/\tau_n \to 0$. Then $\delta/\tau_n$ and $-\delta/\tau_n$ are so close that the optimal limiting power against any alternative sequence $h_n$ with $|h_n| < \delta/\tau_n$ is $\alpha$.

Perhaps the reader can be more easily convinced of these assertions in the $N(\theta, 1)$ model, where explicit expressions for the power of the UMP test exist, but the previous arguments apply to more general models. Thus, the two asymptotic approaches previously considered in this section are in essence the most general.

**Acknowledgment.** Special thanks to Erich Lehmann for some helpful discussion.

DEPARTMENT OF STATISTICS
STANFORD UNIVERSITY
STANFORD, CALIFORNIA 94305-4065
USA
E-MAIL: romano@stat.stanford.edu